\documentclass[11pt]{amsart}

\theoremstyle{definition}
\newtheorem{theorem}{Theorem}[section]

\newtheorem{lemma}[theorem]{Lemma}
 
\newtheorem{prop}[theorem]{Proposition}
\newtheorem{definition}[theorem]{Definition}
\newtheorem{example}[theorem]{Example}

\theoremstyle{remark}
\newtheorem{remark}[theorem]{Remark}

\usepackage{framed}

\usepackage{geometry}
\usepackage{amsfonts}
\usepackage{enumerate}
\usepackage{graphicx}
\usepackage{epstopdf}

\usepackage{float}
\usepackage{algorithmic}
\usepackage{hyperref}
\ifpdf
  \DeclareGraphicsExtensions{.eps,.pdf,.png,.jpg}
\else
  \DeclareGraphicsExtensions{.eps}
\fi
\numberwithin{theorem}{section}
\usepackage{caption}
\usepackage{multirow}
\usepackage{amssymb}
\usepackage{cleveref}
\usepackage{url}
\usepackage{graphicx}
\usepackage{color}
\usepackage{amsopn}
\usepackage{mathtools}

\newcommand{\important}{\emph} 

\newcommand{\DZO}{(\mathbb{C}^*)^n\times(\mathbb{C}^*)^{\theCodim}}
\newcommand{\LagAmbient}{\CC^n\times\PP^\theCodim}

\newcommand{\X}{X}  
\newcommand{\Xo}{{\X^o}}  
\newcommand{\xv}{\mathbf{x}} 

\newcommand{\Z}{Z}  
\newcommand{\Zo}{{\Z^o}}  
\newcommand{\zv}{\mathbf{z}} 

\newcommand{\cH}{\mathcal{H}}  
\newcommand{\Hpoly}{H}  

\newcommand{\pointP}{\mathbf{p}} 
\newcommand{\pointQ}{\mathbf{q}} 

\newcommand{\F}{F} 
\newcommand{\M}{M} 
\newcommand{\J}{J} 

\newcommand{\Fhat}{\hat \F} 
\newcommand{\Mhat}{\hat \M} 
\newcommand{\Jhat}{\hat \J} 

\newcommand{\MLag}{\hat \M_{\Lag}} 

\newcommand{\Jac}{\nabla}

\newcommand{\logOneForm}{\ell_{\pmb{\mu}}(\zv)} 

\newcommand{\theY}{{\texttt{y}}}
\newcommand{\bb}{{b}}

\newcommand{\CC}{\mathbb{C}}
\newcommand{\PP}{\mathbb{P}}
\newcommand{\strata}{S}
\DeclareMathOperator{\Eu}{Eu}
\DeclareMathOperator{\diag}{diag}
\DeclareMathOperator{\Lag}{Lag}
\DeclareMathOperator{\codim}{codim}
\newcommand{\theCodim}{c}  

\DeclareMathOperator{\MLdeg}{MLdegree}

\newcommand{\ThePoint}{{\texttt{ThePoint}}}
\newcommand{\Directory}{{\texttt{Directory}}}
\newcommand{\TheVariety}{{\texttt{TheVariety}}}
\newcommand{\RemovalMLDegree}{{\texttt{RemovalMLDegree}}}
\newcommand{\MLDegreeWitnessCollection}{{\texttt{MLDegreeWitnessCollection}}}
\newcommand{\newMLDegreeWitnessCollection}{{\texttt{newMLDegreeWitnessCollection}}}
\newcommand{\MLDegreeWitnessSet}{{\texttt{MLDegreeWitnessSet}}}
\newcommand{\newMLDegreeWitnessSet}{{\texttt{newMLDegreeWitnessSet}}}
\newcommand{\saveWitnessCollectionConfiguration}{{\texttt{saveWitnessCollectionConfiguration}}}
\newcommand{\getWitnessCollection}{{\texttt{getWitnessCollection}}}
\newcommand{\reclassifyWitnessPoints}{{\texttt{reclassifyWitnessPoints}}}
\newcommand{\homotopyRemovalMLDegree}{{\texttt{homotopyRemovalMLDegree}}}

\newcommand{\temporaryFileName}{{\texttt{temporaryFileName}}}
\newcommand{\mkdir}{{\texttt{mkdir}}}
\newcommand{\WitnessRing}{{\texttt{WitnessRing}}}
\newcommand{\SortPoints}{{\texttt{SortPoints}}}
\newcommand{\WitnessSets}{{\texttt{WitnessSets}}}

\newcommand{\Configuration}{{\texttt{Configuration}}}
\newcommand{\randomValue}{{\texttt{randomValue}}}
\newcommand{\RandomCoefficients}{{\texttt{RandomCoefficients}}}
\newcommand{\saturate}{{\texttt{saturate}}}
\newcommand{\ideal}{{\texttt{ideal}}}

\newcommand{\idealI}{{\texttt{I}}}
\newcommand{\newMLDegreeVariety}{{\texttt{newMLDegreeVariety}}}
\newcommand{\MLDegreeVariety}{{\texttt{MLDegreeVariety}}}
\newcommand{\DefiningEquations}{{\texttt{DefiningEquations}}}
\newcommand{\DataOne}{{\texttt{Data1}}}
\newcommand{\DataZero}{{\texttt{Data0}}}
\newcommand{\Hyperplanes}{{\texttt{Hyperplanes}}}
\newcommand{\MLDegrees}{{\texttt{MLDegrees}}}
\newcommand{\removalMLDegree}{{\texttt{removalMLDegree}}}
\newcommand{\solveRemovalMLDegree}{{\texttt{solveRemovalMLDegree}}}
\newcommand{\mlObstructionFunction}{{\texttt{mlObstructionFunction}}}

\newcommand{\false}{{\texttt{false}}}
\newcommand{\true}{{\texttt{true}}}

\newcommand{\BertiniMTwo}{{\textsc{Bertini.M2}}}
\newcommand{\Bertini}{{\textsc{Bertini}}}
\newcommand{\MacaulayTwo}{{\textsc{Macaulay2}}}

\newcommand{\rk}{{\texttt{k}}}
\newcommand{\mm}{{\texttt{m}}}
\newcommand{\optTo}{{\texttt{=>}}}

\numberwithin{equation}{section} 

\begin{document}

\title{
Implementations of symbolic-numeric algorithms computing Euler obstruction functions using
maximum likelihood degrees}
\title{Solving the likelihood equations to compute Euler obstruction functions}

\author[J.~I.~Rodriguez]{Jose Israel Rodriguez}
\address{Department of Statistics\\
         University of Chicago\\
         Chicago, IL 60637\\         
         USA}
\email{JoIsRo@uchicago.edu}
\urladdr{{http://home.uchicago.edu/\~joisro}}

\maketitle
\begin{abstract}
Macpherson defined Chern-Schwartz-Macpherson (CSM) classes by introducing the (local) Euler obstruction function, which is an integer valued function on the variety that is constant on each stratum of a Whitney stratification of an algebraic variety.
 By understanding the Euler obstruction function, one gains insights about a singular algebraic variety. 
It was recently shown by the author and B. Wang, how to compute these functions using maximum likelihood degrees. 
 This paper discusses a symbolic and a numerical implementation of algorithms to compute the Euler obstruction at a point.
 Macaulay2 and Bertini are used in the implementations. 
\end{abstract}




\section{Introduction }
Studying singularities of algebraic varieties is of great interest in applied and computational algebraic geometry. For example, in applications the singular locus is important when finding the closest point to an algebraic variety. In computational algebraic geometry, the singular locus can lead to bottlenecks of an algorithm. 
One way to understand a singular algebraic variety is by stratifying it into locally closed sets called Whitney strata. Then, for each stratum one considers the local information at a point.

To study how the closures of these strata interact with one another, the Euler obstruction function defined in \cite{Mac74}  by Macpherson is considered. 
 Informally, this function gives a measure of the singularity of a stratum. For an equivalent definition using Euler characteristics of complex links see \cite{BLS00,Dim04}.  
In \cite{RW2017b},
 it is shown that the Euler obstruction function at a point is given by an alternating sum of maximum likelihood degrees, which will be reviewed in the next section. 
It is in this framework, where we will develop our algorithms.
The first sections of this paper recall the definition of ML degree, likelihood equations, and removal ML degrees to provide the statement of Theorem~\ref{theorem:result} \cite{RW2017b}, which enables us to compute the Euler obstruction function. 
In Section \ref{sec:Example}, a concrete example is provided to set notation and motivate the work. 
 The main result of this paper is a description of a Macaulay2 \cite{M2} package in Section \ref{sec:package} that has implemented the algorithms of  \cite{RW2017b}.
The goal is to provide the research community tools to compute interesting examples to drive new areas of study. 
In turn, this drives the search for new algorithms.\footnote{The author is thankful for the helpful comments and suggestions of Botong Wang and Xiping Zhang
 at ``Singularities in the Midwest, V''.
}
The package is available at:
\begin{quote}
\url{https://github.com/JoseMath/MaximumLikelihoodObstructionFunction}
\end{quote}


\section{Maximum likelihood degrees}\label{sec:mlDegrees}
\subsection{Maximum likelihood degrees of very affine varieties}\label{ss:veryAffine}
Consider an irreducible affine variety $\Z$ of $\CC^n$.
The set of points of $\Z$ with nonzero coordinates 
is denoted by $\Zo$ and said to be 
the \important{underlying very affine variety of} $\Z$.
The underlying very affine variety of $\Z$ is a subvariety of $(\CC^{*})^n$.
Consider the logarithmic $1$-form on $(\CC^*)^n$ given by 
\[
\logOneForm:=\mu_1\log{z_1}+\mu_2\log{z_2}+\cdots+\mu_n\log{z_n},
\]
where $\pmb\mu=(\mu_1,\dots,\mu_n)\in \CC^n$. 
The gradient of this one form is 
\[
\Jac\logOneForm\coloneqq
\begin{bmatrix}
\mu_1/z_1& \dots &\mu_n/z_n
\end{bmatrix}.
\]
 Let $\zv$ denote a point of $\Zo$.
If $\zv$ is a regular point of $\Zo$,
then the one form restricted to $\Zo$ is said to have a  \important{critical point}  $\zv$ 
if  $\Jac\logOneForm$ at $\zv$ is orthogonal to the tangent space of $\Zo$ at $\zv$. 

 \begin{definition}\label{definition:mlDegree}
The  \important{maximum likelihood degree} (ML degree) of the very affine variety $\Zo$,
is defined 
to be the number of critical points of $\logOneForm$ on $\Zo$
 for general $\mu_i$ 
 and is 
 denoted by $\MLdeg(\Zo)$.
The maximum likelihood degree of an affine variety $\Z$
is defined to be the ML degree of the underlying very affine variety $\Zo$.
\end{definition}

The notion of ML degree was first introduced in \cite{CHKS06,HKS05}. 
The name  ``maximum likelihood" comes from statistics, where the log likelihood function is the 1-form $\logOneForm$, with $\mu_i$ denoting the number of times event $i$ is observed. 
For a more geometric interpretation Definition~\ref{definition:mlDegree} see \cite{Huh13}, and for a survey of results see \cite{HS2014}.
Moreover,   the \emph{Gaussian degree} \cite{FK00} is in some cases equivalent to the ML degree, and 
the data singular locus for maximum likelihood estimation is studied in \cite{EmilJose2017}.
In addition, the maximum likelihood degree appears in other contexts, including Gaussian graphical models \cite{Uhler12}, variance component models \cite{GDP12}, and in missing data \cite{HS09}. 

Our convention is that the ML degree of an empty set is zero. 

\subsection{Likelihood equations}
The critical points of $\logOneForm$ are defined by a zero dimensional variety and can be found by solving a system of polynomial equations. 
We determine the maximum likelihood degree by determining the degree of this zero dimensional variety. 
Let $\F\subset\CC[z_1^{\pm},\dots,z_n^{\pm}]$ denote a set of generators of the ideal of $\Zo$,
$\J$ denote the ideal of the singular locus of $\Zo$, and
 $\M$ denote the ideal of the $1+\codim(\Zo)$ minors of 
\begin{equation}\label{eq:M}
\begin{bmatrix}
\nabla \logOneForm\\
\Jac \F\end{bmatrix},
\end{equation}
where $\Jac \F$ is the matrix of partial derivatives of $\F$.
The variety of  $\saturate(\M+\F,\J)$
 is the set of critical points. 

\begin{definition}\label{definition:LE}
The \important{likelihood equations} of $\Zo$ (with respect to $\logOneForm$)
are defined by setting each element of a set of generators of the ideal 
 $\saturate(\M+\F,\J)$
to be zero. 
\end{definition}

\begin{prop}
For generic $\pmb\mu$, 
the degree of the variety of the ideal of the likelihood equations is the maximum likelihood degree of $\Zo$.
\end{prop}

Often, it is easier to work with the affine variety $\Z$ rather than $\Zo$. 
Let $\Fhat$ denote a set of generators of the ideal of $\Z$,
$\Jhat$ denote a set of generators of the ideal of the singular locus of $\Z$, 
and 
$\Mhat$ denote a set of generators of the $(\codim\Zo+1)$-minors of 
\begin{equation}\label{eq:Mhat}\begin{bmatrix}
\nabla \logOneForm\\ 
\Jac \Fhat\end{bmatrix}\diag(\begin{bmatrix}z_1&\dots &z_n\end{bmatrix}).
\end{equation}
 
\begin{lemma}
Setting each element of a set of generators of the ideal 
 \[
 \saturate(\ideal\, \Mhat+\ideal\, \Fhat,\ideal\, \Jhat*\ideal\, (z_1z_2\cdots z_n))
 \] 
to zero defines likelihood equations of $\Zo$. 
\end{lemma}
\begin{remark}
The definition of the ML degree of a projective variety in \cite{HKS05} agrees with the definition of the ML degree of a very affine variety when we restrict the projective variety to the affine chart where the coordinates sum to one. 
\end{remark}

The likelihood equations are often an overdetermined system of equations. 
For numerical computation we prefer for the number of equations to equal the dimension of the ambient space. 
Suppose $\Fhat$ is a set of $\theCodim\coloneqq\codim\Zo$  generators of the ideal of the ideal $\Z$.
Let $\{\lambda_0,\dots,\lambda_{\theCodim}\}$ denote a set of indeterminants called \important{Lagrange multipliers}.
Let $\MLag$ denote the following set of $n$ polynomials
\begin{equation}\label{eq:MLag}
\begin{bmatrix}\lambda_0&\lambda_1&\dots &\lambda_{\theCodim}\end{bmatrix}
\begin{bmatrix} 
\nabla \logOneForm\\ 
\Jac \F\end{bmatrix}
\diag(\begin{bmatrix}z_1&\dots &z_n\end{bmatrix}).
\end{equation}

\begin{definition}\label{definition:lagEq} 
The \important{Lagrange likelihood equations} are a system of equations 
defined by $\Fhat=0$ and $\MLag=0$.
\end{definition}

The variety of  $\ideal\,\Fhat+\ideal\,\MLag$ is in the product space $\LagAmbient$.
Let $\DZO$ be the dense Zariski open set of $\LagAmbient$
defined by $\lambda_0=1,\prod_{i=1}^nz_i\prod_{j=1}^{\theCodim}\lambda_j\neq0$.

\begin{prop}
For generic $\pmb\mu$, 
the projection 
to $(\mathbb{C}^*)^n$
of $\DZO$  intersected with  the variety of the Lagrange likelihood equations  
is the variety of the likelihood equations of $\Zo$.
The degree of the projection is one. 
\end{prop}
These Lagrange likelihood equations were studied in \cite{GR13}.

\subsection{Computing removal ML degrees and Euler obstructions}\label{ss:removalMLDegree}
Let $\Xo$ denote a very affine variety of $(\CC^{*})^n$, and let $\pointP$ denote a point in $(\CC^{*})^n$.
Let $\cH_1,\dots,\cH_n$ denote general hyperplanes containing the point $\pointP$, and let 
$\Hpoly_1(\xv),\dots,\Hpoly_n(\xv)$ denote affine linear polynomials defining these hyperplanes. 

We denote embedding of $\Xo$ to $(\CC^{*})^{n+1}$ via $y=\Hpoly_1(\xv)$ by $\Xo\setminus \cH_1$.
\begin{definition}\label{definition:removalMLDegree}
For $k=0$,  let $\Zo=\Xo$, otherwise let $\Zo=(\Xo\setminus \cH_1)\cap\left(\cap_{i=2}^k\cH_i\right)$.
The $k$-th removal ML degree of $\Xo$ with respect to the point $\pointP$ 
is defined to be the ML degree of $\Zo$ and is denoted by $r_k(\Xo,\pointP)$.
\end{definition}

With removal ML degrees, we can compute the Euler obstruction by the following theorem. 

\begin{theorem}[\cite{RW2017b}]\label{theorem:result}
The signed alternating sum of removal ML degrees with respect to the point $\pointP$ in $\Xo$ of dimension $d$
equals to value of the Euler obstruction function at a point $\pointP$, i.e.,
\[
(-1)^{d}\Eu_\Xo(\pointP)=\sum_{k=0}^{d+1} (-1)^kr_k(\pointP,\Xo).
\]
\end{theorem}

With the symbolic implementation, we compute the removal ML degrees by solving the likelihood equations using Grobner basis.
With the numerical implementation, we will use \emph{coefficient parameter homotopies} \cite[Theorem 7.1.1]{SW05}.
For fixed $k$, we choose a generic $k\times n$ matrix $[\gamma_{ij}]_{k\times n}$, 
and let $\mathcal{L}_\bb$ 
denote the following family of linear spaces defined by
\begin{align}\label{eq:familyLinearSpaces}
\begin{bmatrix}
\Hpoly_1(\xv) \\ \vdots \\ \Hpoly_k(\xv)
\end{bmatrix} &=\begin{bmatrix}\gamma_{ij}\end{bmatrix}_{k\times n}\begin{bmatrix}
x_1\\
 \vdots\\
 x_n 
\end{bmatrix}
-\begin{bmatrix}
b_1\\
 \vdots\\
 b_k
\end{bmatrix}.
\end{align}

Let $\Z=\X\cap\mathcal{L}_\bb$. 
For general $[\gamma_{ij}]_{k\times n}$ in \eqref{eq:familyLinearSpaces}, 
 the Lagrange likelihood equations of $\Z$ define a coefficient parameter homotopy with \emph{parameters} $\bb=(b_1,\dots,b_k)$,
Lagrange multipliers $\lambda_0,\dots,\lambda_\theCodim$, and \emph{primal variables} $x_1,\dots,x_n,y$ for $k>0$.

\begin{theorem}[Corollary 4.6 \cite{RW2017b}]\label{theorem:homotopy}
If $\gamma_{ij}$ are general and $\Z=\X\cap\mathcal{L}_\bb$, 
then for  a parameter homotopy with target parameters $([\gamma_{ij}]_{k\times n}[\pointP])$,
the number of regular endpoints not in the coordinate hyperplanes 
equals the $k$-th removal ML degree of $\X$ with respect to the point 
$\pointP$.
\end{theorem}

In summary, we solve the Lagrange likelihood equations for a general choice of parameters, thereby determining the removal ML degrees with respect to a general point. 
Then, we use a parameter homotopy to determine the removal ML degrees for any other point $\pointP$ of interest.

\section{Whitney Sombrilla}\label{sec:Example}

The illustrative example in this section comes from a Whitney umbrella in $\CC^3$ defined by 
$x_1^2-x_2^2x_3=0$ that has been translated by $(1,1,1)$.
We call this translation the \important{Whitney sombrilla} and denote it by $X$.
The defining equation is $f\coloneqq(x_1-1)^2-(x_2-1)^2(x_3-1)=0$. 
Note that the original Whitney umbrella and the Whitney sombrilla have different underlying very affine varieties. 
Indeed, the underlying very affine variety of the Whitney umbrella is smooth while the underlying very affine variety of the Whitney sombrilla is not. 
A Whitney stratification of $\X$ is given by the regular points; the singular points with $\{(1,1,1)\}$ removed; and $\{(1,1,1)\}$.
We denote these strata by $\strata_1,\strata_2,\strata_3$ respectively. Let $\strata_0$ denote $\CC^3\setminus \X$.
Let $\pointP_i$ be a point in $\strata_i\cap(\CC^*)^3$.  
The removal ML degrees of these points and their respective Euler obstructions are below.
Note that the defining equations of $\X\setminus\cH_{1},\,
\left(\X\setminus\cH_{1}\right)\cap\cH_{2},\,
\left(\X\setminus\cH_{1}\right)\cap(\cH_2\cap\cH_3)$ 
 are given by $f=0$ and
  $  \{H_1=y\},
 \{
 H_1=y, H_2=0 \},
 \{
 H_1=y,
 H_2=0,
 H_3=0\}$ respectively. 
\[
\begin{tabular}{r||c|c|c|c||c}
 & $k=0$ & $k=1$ & $k=2$ &$k=3$ &$\Eu$ \\
\hline\hline
$\pointP_0=(3,2,1)$ & 3 & 10 &  10&	3 & 0\tabularnewline
\hline
$\pointP_1=(3,3,2)$ & 3 & 10 &  10&	2& 1\tabularnewline
\hline
$\pointP_2=(1,1,2)$ & 3  &10 &  10&	1& 2\tabularnewline
\hline
$\pointP_3=(1,1,1)$ & 3  &10 &  9 &       1& 1\tabularnewline
\hline
\end{tabular}
\]

\section{Using the package}\label{sec:package}
The package computes the Euler obstruction function of a very affine
variety $\Xo\subset(\CC^*)^n$ at the point $\pointP$.
This is done by computing the $k$-th removal ML degrees of $\Xo$ with respect to the point $\pointP$ for $k=0,1,\dots,\dim \Xo+1$.
We define a new type of mutable hash table
called the \RemovalMLDegree, which is used to store the results of the computations of the removal ML degrees. 

The package takes two approaches to computing removal ML degrees. 
The first approach uses symbolic computation similar to the approaches in the foundational paper \cite{HKS05}.
The second approach uses homotopy continuation \cite{allgower1993continuation}.

In each case, our algorithms are probabilistic and there exists a open Zariski dense set such that choices of random values will produce the true answer. 
How we generate random values can be changed in the {\Configuration} when loading the package (see help {\randomValue}). 
The default has the configuration of {\RandomCoefficients} set to $(1,30102)$ so that the random values are produced by 
\texttt{random(1,30102)}.

\subsection{Symbolic computation}\label{ss:symbolic}

\subsubsection{Preprocess}
We assume the variety $\Xo$ is defined by an ideal $\idealI$, but we store information about this variety in a new type of mutable hash table called {\MLDegreeVariety}. 
To create this mutable hashtable 
we use the method {\newMLDegreeVariety}. 
See Table~\ref{table:MLDegreeVariety} 
for a a description of the keys of the hashtable. 

\begin{leftbar}
\begin{verbatim}
i1 :  loadPackage"MaximumlikelihoodObstructionFunction";
i2 :  R=QQ[p1,p2,p3];
i3 :  I=ideal((p-1)1^2-(p2-1)^2*(p3-1));
i4 : 	L=newMLDegreeVariety(I)	
o4 = MLDegreeVariety{...4...}
o4 : MLDegreeVariety
\end{verbatim}
\end{leftbar}

\subsubsection{Solving}\label{ss:solvingSym}
Let $\pointP$ denote a point in $(\CC^*)^n$.
If the point $\pointP$ is not  specified by the user, then
$\pointP$ is set to $(1,\dots,1)$. 
For $k=0,1,\dots,\dim \Xo+1$, the package computes the $k$-th removal ML degree of $\Xo$ with respect to $\pointP$.
To store this information we introduce a type of mutable hash table
called $\RemovalMLDegree$.
This hashtable has three important keys $\MLDegrees$, $\ThePoint$, and $\TheVariety$ described in Table \ref{table:RemovalMLDegree}.
              
The method $\solveRemovalMLDegree$,
solves the likelihood equations (Definition~\ref{definition:LE}) 
and stores the degree in $\RemovalMLDegree$ by appending $\rk\optTo\mm$ to the list under the key $\MLDegrees$, where $\rk$ indexes the values of the removal ML degree $\mm$.
This is is the most difficult step in the computation. 
\begin{leftbar}
\begin{verbatim}
i5 : 	P={1,1,1};
i6 : 	M=newRemovalMLDegree(L,P)
o6 = RemovalMLDegree{...3...}
o6 : RemovalMLDegree
i7 : 	solveRemovalMLDegree M
o7 = {3, 10, 9, 1}
\end{verbatim}
\end{leftbar}

\subsubsection{Extracting information}
Once all of the removal ML degrees are computed we can extract the information using the methods $\removalMLDegree$ or $\mlObstructionFunction$.
The former return lists of computed removal ML degrees where the $k$-th element of the list is the $k$-th removal ML degree of $\Xo$ at $\pointP$. 
The $\mlObstructionFunction$ returns the alternating sum of computed ML degrees.

\begin{leftbar}
\begin{verbatim}
i8 :  removalMLDegree(M)
o8 = {3, 10, 9, 1}
i9 : 	mlObstructionFunction M
o9 = 1
\end{verbatim}
\end{leftbar}

\subsection{Numeric computation}\label{ss:numeric}
In this subsection we compute removal ML degrees using homotopy continuation with the numerical algebraic geometry software $\Bertini$ \cite{Bertini}.
We use methods of the $\MacaulayTwo$ package $\BertiniMTwo$ \cite{Bertini4M2} to manipulate the input files. 

\subsubsection{Preprocess}
Let $\pointQ$ denote a general point in $(\mathbb{C}^*)^n$. 
In this step we compute the $k$-th removal ML degrees of $\Xo$ at $\pointQ$ by solving the equations in Theorem~\ref{theorem:homotopy} for a general choice of parameters. 

We assume the variety $\Xo$ is given to us by an ideal $\idealI$, but we store information in a new type of mutable hash table
called {\MLDegreeWitnessCollection}. 
To create this mutable hashtable 
we use the method {\newMLDegreeWitnessCollection}. 
The keys of the resulting hashtable  include those in 
Table~\ref{table:MLDegreeVariety} and Table~\ref{table:MLDegreeWitnessCollection}.

\newcommand{\Degree}{\texttt{Degree}}
\newcommand{\WitnessPoints}{\texttt{WitnessPoints}}
\newcommand{\VariableGroups}{\texttt{VariableGroups}}

\begin{leftbar}
\begin{verbatim}
i10 : loadPackage("MaximumlikelihoodObstructionFunction",Reload=>true,
Configuration=>{"RandomCoefficients"=>CC})
o10 = MaximumlikelihoodObstructionFunction
o10 : Package
i11 : 	s = temporaryFileName() | "/";
i12 : 	mkdir s;
i13 : 	R=CC[p1,p2,p3];
i14 : 	I=ideal((p-1)1^2-(p2-1)^2*(p3-1))
i15 : 	WC=newMLDegreeWitnessCollection(I,d,s)
o15 = MLDegreeWitnessCollection{...10...}
o15 : MLDegreeWitnessCollection
\end{verbatim}
\end{leftbar}

\subsubsection{Computing witness sets}\label{ss:ComputingWS}
We solve the Lagrange likelihood equations (Definition~\ref{definition:lagEq}) for removal degrees of $\Xo$ with respect to a generic point using the method $\newMLDegreeWitnessSet$. 
This computes the removal ML degrees for $\X$ at a generic point. 
We store this information in a directory determined by the key $\Directory$ of  $\MLDegreeWitnessCollection$. 
The most intensive part of the computation is this step were we compute witness sets. 
To avoid repeating this step we are able to load a witness collection from a saved directory using the method 
$\getWitnessCollection$ if the configuration information of the witness collection is saved using 
$\saveWitnessCollectionConfiguration$.
The $\MLDegreeWitnessSet$ has its keys and values described in Table~\ref{table:MLDegreeWitnessSet}.
\begin{leftbar}
\begin{verbatim}
i16 : 	newMLDegreeWitnessSet(WC)
o16 = {3, 10, 10, 3}
i17 :  saveWitnessCollectionConfiguration(WC,s)
\end{verbatim}
\end{leftbar}

\subsubsection{Solving}\label{ss:solvingNum}
With the method $\homotopyRemovalMLDegree$, 
we use a parameter homotopy from Theorem~\ref{theorem:homotopy} 
to determine the Euler obstruction at a point. 
This step is much faster than the previous step.

\begin{leftbar}
\begin{verbatim}
i18 :  P={1,1,1}
i19 :  M=newRemovalMLDegree(WC,P)
o19 = RemovalMLDegree{...4...}
o19 : RemovalMLDegree
i20 : 	homotopyRemovalMLDegree M
o20 = {3, 10, 9, 1}
\end{verbatim}
\end{leftbar}

\subsubsection{Extracting information}
Since we are working with floating point arithmetic, one must take care when classifying points in the coordinate hyperplanes. The method $\reclassifyWitnessPoints$ allows us to change the tolerances as we like. 

\begin{leftbar}
\begin{verbatim}
--Tolerance is too tight compared to the precision used.
i21 : 	reclassifyWitnessPoints(M,1e-300)                                                                              
o21 = {3, 10, 10, 2} --Incorrect answer
--This tolerance yields the correct answer.
i22 : 	reclassifyWitnessPoints(M,1e-6)                                                       
o22 = {3, 10, 9, 1}  
\end{verbatim}
\end{leftbar}

\section{Motivating examples}
Matrices with rank constraints give a wide collection of interesting examples and have been studied in \cite{HRS14,KRS15}. Recent work in \cite{RW2017b}, uses Euler obstructions to prove results about ML degrees that were motivated by previous computations, and these obstructions have also been studied in \cite{GGR2017,Zhang2017}.
The first example considers $\X_1\subset\CC^4$ defined by the determinant of the matrix on the left, and the second example considers $\X_2\subset\CC^5$ defined by the determinant of the matrix on the right \eqref{eq:m}. 
We let $\strata_0$ denote the ambient space with $\X_i$ removed; 
$\strata_1$ denote the regular points of $\X_i$;
$\strata_2$ denote the singular points in $\X_i$. 
For the benchmarks below: 
the $k$th column records  $r_k(X_i,\pointP)$ for $\pointP=\pointP_0,\pointP_1,\pointP_2$ and the wall time 
of the computing $k$th witness set step in Section~\ref{ss:ComputingWS}. 
The $\Eu$-column records the Euler obstruction at $\pointP_i$.  The Time-Sym and Time-Num columns record
 the wall time of the solving step in Section~\ref{ss:solvingSym} and  Section~\ref{ss:solvingNum} respectively.

\begin{equation}\label{eq:m}
\begin{bmatrix}
x_1 & x_2 &x_3\\
x_2 & x_3 &x_4\\
x_3 &x_4 &1
\end{bmatrix}
\quad \quad\quad
\begin{bmatrix}
x_1 & x_2 &x_3\\
x_2 & x_4 &x_5\\
x_3 &x_5 &1
\end{bmatrix}
\end{equation}
\[
\begin{tabular}{c||c|c|c|c|c||c|c|c}
$\X_1$ & $k=0$ & $k=1$ & $k=2$ &$k=3$&$k=4$ &$\Eu$& Time-Sym  &Time-Num\\
\hline\hline
$\pointP_0=(7,5,3,2)$ & 0 & 16 &  31&	18 & 3 &0 & 58s& 1s\tabularnewline
\hline
$\pointP_1=(2,1,1,1)$ & 0 & 16 &  31&	18 & 2 &1 & 56s& 1s\tabularnewline
\hline
$\pointP_2=(1,1,1,1)$ & 0 & 16 &  31&	16 & 1 &0 & 62s& 1s\tabularnewline
\hline
\hline
Time-Compute WS& 20s& 46s& 52s& 27s& 3s\tabularnewline
\end{tabular}
\]
\[
\begin{tabular}{c||c|c|c|c|c|c||c|c|c}
$\X_2$ & $k=0$ & $k=1$ & $k=2$ &$k=3$&$k=4$&$k=5$ &$\Eu$& Time-Sym &Time-Num\\
\hline\hline
$\pointP_0=(1,2,3,5,7)$ & 0 & 16 &  47&	49 & 21 &3& 0 & 1620s& 1s\tabularnewline
\hline
$\pointP_1=(1,1,1,1,2)$& 0 & 16 &  47&	49 & 21 &2& 1 &1991s& 1s\tabularnewline
\hline
$\pointP_2=(1,1,1,1,1)$& 0 & 16 &  47&	49 & 19 &1& 0& 2002s& 1s\tabularnewline
\hline
\hline
Time-Compute WS & 8s& 197s& 321s& 90s& 68s& 9s\tabularnewline
\end{tabular}
\]

\bibliographystyle{abbrv}
\bibliography{BIB_ComputeEO}


\begin{table}[hbt!]
\centering
\caption{Keys of $\MLDegreeVariety$}
\label{table:MLDegreeVariety}
\[
\begin{tabular}{|p{3.5cm}|p{12cm}|}
\hline 
  Key and Value &  Details \\
\hline\hline
$\DefiningEquations$ \newline Ideal
&  After saturating by the coordinate hyperplanes, this variety defines $\X$.
\tabularnewline
\hline
$\DataZero$, $\DataOne$ \newline List
& 
Random values taking the role of $\pmb\mu$ in the likelihood equations.
For $k=0$,  $\DataZero$  is used, otherwise $\DataOne$ is used. 
\tabularnewline
\hline
$\Hyperplanes$ \newline List 
&
The $i$-th element defines the gradient of the polynomial $H_i(x)$, i.e., 
 the normal vector to the hyperplane $\cH_i$. 
\tabularnewline
\hline
\end{tabular}
\]
\end{table}

\begin{table}[hbt!]
\centering
\caption{Keys of $\RemovalMLDegree$}
\label{table:RemovalMLDegree}
\[
\begin{tabular}{|p{3.5cm}|p{12cm}|}
\hline
%
$\MLDegrees$ \newline List
&   List of options $\rk\optTo \mm$ where $\mm$ is the value of the $\rk$-th removal ML degree. 
This is initialized to be the empty set and is determined by using the method $\solveRemovalMLDegree$.
\tabularnewline
\hline
$\ThePoint$ \newline List
& 
$\pointP$ with the default set to $(1,1,\dots,1)$. 
\tabularnewline
\hline
$\TheVariety$ \newline \MLDegreeVariety
& 
Defines the variety of interest.  
\tabularnewline
\hline
\end{tabular}
\]
\end{table}

\begin{table}[hbt!]
\centering
\caption{Keys of $\MLDegreeWitnessCollection$ (see also Table~\ref{table:MLDegreeVariety})}
\label{table:MLDegreeWitnessCollection}
\[
\begin{tabular}{|p{3.5cm}|p{12cm}|}
\hline
%
$\Directory$ \newline 
String 
& The string is given by $\temporaryFileName()$. The directory is created if it it doesn't exist using $\mkdir$ and stores the output of 
$\Bertini$ runs. 
\tabularnewline
\hline
$\WitnessRing$ \newline Ring 
& This ring has three variable groups: the coordinates of $\Xo$ and $\theY$,  the Lagrange multipliers, and the parameters.
\tabularnewline
\hline
$\SortPoints$ \newline Function
& The function takes a list of coordinates of a point as an input  and returns $\true$ if the point is in $\DZO$ and $\false$ otherwise. 
\tabularnewline
\hline 
$\WitnessSets$ \newline List of options
&This key is initialized to be the empty set. 
The method $\newMLDegreeWitnessSet$ is used to create a new type of mutable hash table 
called a $\MLDegreeWitnessSet$ (Table~\ref{table:MLDegreeWitnessCollection}).
\tabularnewline
\hline 
\end{tabular}
\]
\end{table}

\begin{table}[hbt!]
\centering
\caption{Keys of $\MLDegreeWitnessSet$}
\label{table:MLDegreeWitnessSet}
\[
\begin{tabular}{|p{3.5cm}|p{12cm}|}
\hline
%
$\WitnessPoints$ \newline Sequence of \newline 2 elements
& The first element is a list of list of coordinates of points.  The second element is a list of booleans, that are the result of applying $\SortPoints$ to the first element. 
\tabularnewline
\hline
$\Degree$ \newline Integer 
& The number of points outside of the coordinate hyperplanes. 
This is determined by counting the number of times true appears in the second element of \WitnessPoints.
\tabularnewline
\hline
$\DefiningEquations$ \newline List of 3 elements
& The first element are defining equations of $\X$, the second are defining equations for $\cH_1,\dots,\cH_k$, and the third are the equations \eqref{eq:MLag}.
\tabularnewline
\hline 
$\VariableGroups$ \newline List of 3 elements
& The first element is the set of Lagrange multipliers and is a homogeneous variable group; 
the second element includes $x_1,\dots,x_n$ and $y$ if $k\neq 0$; the third is the set of parameters $\bb$ used 
Theorem~\ref{theorem:homotopy}.
\tabularnewline
\hline 
\end{tabular}
\]
\end{table}

\end{document}